\documentclass[journal]{IEEEtran}

\usepackage{myStyleIEEE}
\makenomenclature

\newenvironment{ldescription}[1]
  {\begin{list}{}%
   {\renewcommand\makelabel[1]{##1\hfill}%
   \settowidth\labelwidth{\makelabel{#1}}%
   \setlength\leftmargin{\labelwidth}
   \addtolength\leftmargin{\labelsep}}}
  {\end{list}}

\IEEEoverridecommandlockouts

\begin{document}

\title{Stochastic Unit Commitment in Low-Inertia Grids}

\renewcommand{\theenumi}{\alph{enumi}}

\newcommand{\uros}[1]{\textcolor{magenta}{$\xrightarrow[]{\text{Uros}}$ ``#1''}}
\newcommand{\vaggelis}[1]{\textcolor{blue}{$\xrightarrow[]{\text{Vaggelis}}$ ``#1''}}
\newcommand{\stefanos}[1]{\textcolor{red}{$\xrightarrow[]{\text{Stefanos}}$ ``#1''}}
\newcommand{\matthieu}[1]{\textcolor{pPurple}{$\xrightarrow[]{\text{Matt}}$ ``#1''}}

\author{Matthieu~Paturet, Uros~Markovic,~\IEEEmembership{Student~Member,~IEEE,}
        Stefanos Delikaraoglou,~\IEEEmembership{Member,~IEEE,} \\
        Evangelos Vrettos,~\IEEEmembership{Member,~IEEE,}
        Petros Aristidou,~\IEEEmembership{Member,~IEEE,}
        and~Gabriela~Hug,~\IEEEmembership{Senior~Member,~IEEE}
        }

\maketitle
\IEEEpeerreviewmaketitle

\begin{abstract}
In this paper, the Unit Commitment (UC) problem in a power network with low levels of rotational inertia is studied. Frequency-related constraints, namely the limitation on Rate-of-Change-of-Frequency (RoCoF), frequency nadir and steady-state frequency error, are derived from a uniform system frequency response model and included into a stochastic UC that accounts for wind power and equipment contingency uncertainties using a scenario-tree approach. In contrast to the linear RoCoF and steady-state frequency error constraints, the nadir constraint is highly nonlinear. To preserve the mixed-integer linear formulation of the stochastic UC model, we propose a computationally efficient approach that allows to recast the nadir constraint by introducing appropriate  bounds on relevant decision variables of the UC model. For medium-sized networks, this method is shown to be computationally more efficient than a piece-wise linearization method adapted from the literature. Simulation results for a modified IEEE RTS-96 system revealed that the inclusion of inertia-related constraints significantly influences the UC decisions and increases total costs, as more synchronous machines are forced to be online to provide inertial response.

\end{abstract}



\section*{Nomenclature}
The main notation used in this paper is introduced below. Additional symbols are defined in the paper where needed. All symbols are augmented by index $t$ when referring to different time periods.

\vspace{-0cm}
\subsection{Sets and Indices}
\begin{ldescription}{$xxxxx$}
\item [$\ell \in \mathcal{L}$] Set of transmission lines.
\item [$\xi \in \mathcal{E}$] Set of scenarios $\xi = \{c, \omega\}$ including generation outages ($c$) and wind power uncertainty ($\omega$).
\item [$i \in \mathcal{I}$] Set of conventional generation units.
\item [$j \in \mathcal{J}$] Set of converter-based (i.e., wind) generation units.
\item [$n \in \mathcal{N}$] Set of nodes.
\item [$\mathcal{I}_n$] Set of conventional generation units located at bus $n$.
\item [$\mathcal{J}_n$] Set of converter-based units located at bus $n$.
\end{ldescription}

\vspace{-0cm}
\subsection{Decision variables}
\begin{ldescription}{$xxxx$}
\item [$\hat{\delta}_{n}$] Day-ahead voltage angle at node $n$ $[\mathrm{rad}]$.
\item [$\tilde{\delta}_{n\xi}$] Real-time voltage angle at node $n$ in scenario $\xi$ $[\mathrm{rad}]$.
\item [$F_{\xi t}$] Global fraction of total power generated by high-pressure turbines in scenario $\xi$ $[\mathrm{p.u.}]$.
\item [$k_{i\xi t}$] Scaled power gain factor of conventional unit $i$ $\xi$ $[\mathrm{MW}]$.
\item [$l_{n\xi t}^\mathrm{shed}$] Shedding of load at node $n$ in scenario $\xi$ $[\mathrm{MW}]$.
\item [$M_{\xi t}$] Global system inertia in scenario $\xi$ $[\mathrm{p.u.}]$.
\item [$p_{it}$] Day-ahead dispatch of conventional unit $i$ $[\mathrm{MW}]$.
\item [$r_{i\xi t}^{+/-}$] Up-/Downward reserve deployment of unit $i$ in scenario $\xi$ $[\mathrm{MW}]$.
\item [$R_{\xi t}$] Global system droop factor in scenario $\xi$ $[\mathrm{p.u.}]$.
\item [$u_{it}$] Commitment variable of conventional unit $i$.
\item [$w_{jt}$] Day-ahead dispatch of wind power unit $j$ $[\mathrm{MW}]$.
\item [$w_{j\xi t}^\mathrm{spill}$] Wind spillage of unit $j$ in scenario $\xi$ $[\mathrm{MW}]$.
\item [$y_{it}$] Start-up variable of conventional unit $i$ .
\item [$z_{it}$] Shut-down variable of conventional unit $i$ .
\end{ldescription}

\vspace{-0cm}
\subsection{Parameters}
\begin{ldescription}{$xxxxx$}
\item [$\alpha_{i\xi t}$] Outage parameter of conventional unit $i$ in scenario $\xi$.
\item [$\Delta P_{\xi t}$] Size of power outage in scenario $\xi$ $[\mathrm{p.u.}]$.
\item [$\pi_{\xi}$] Probability of occurrence of scenario $\xi$.
\item [$B_{nm}$] Susceptance of transmission line $(n,m)$ $[\mathrm{S}]$.
\item [$C_{i}$] Day-ahead price offer of unit $i$ $[\mathrm{\$/MWh}]$.
\item [$C_{i}^{\mathrm{SU/SD}}$] Start-up/Shut-down price offer of unit $i$ $[\mathrm{\$}]$.
\item [$C_{i}^{\mathrm{+/-}}$] Up-/Downward reserve price offer of unit $i$ $[\mathrm{\$/MWh}]$.
\item [$C^\mathrm{sh}$] Value of lost load $[\mathrm{\$/MWh}]$.
\item [$D_{nt}$] Demand at node $n$ $[\mathrm{MW}]$.
\item [$\overline{f}_{nm}$] Capacity of transmission line $(n,m)$ $[\mathrm{MW}]$.
\item [$\overline{P_{i}} / \underline{P_{i}}$] Active power limits of conventional unit $i$ $[\mathrm{MW}]$.
\item [$R_{i}^\mathrm{+/-}$] Up-/Downward reserve capacity of unit $i$ $[\mathrm{MW}]$.
\item [$R^{\mathrm{U}/\mathrm{D} }_i$] Ramp up/down limits of conventional unit $i$ $[\mathrm{MW/h}]$.
\item [$W_{j\xi t}^*$] Wind power realization of unit $j$ in scenario $\xi$ $[\mathrm{MW}]$.
\item [$W_{d}$] Total capacity of droop control units $[\mathrm{MW}]$.
\item [$W_{v}$] Total capacity of VSM units $[\mathrm{MW}]$.
\end{ldescription}
\vspace{0.3cm}

\section{Introduction} \label{sec:Intro}

With increasing penetration of renewable energy sources, system operators face new challenges in order to ensure power grid stability. One of these challenges is frequency stability due to a loss of generation or a large variation of load. In traditional power systems, synchronous generators (e.g., hydro or steam turbines) provide rotational inertia through stored kinetic energy in their rotating mass (turbine system and rotor). This energy is important to stabilize the system as it ensures slower frequency dynamics and reduces the Rate of Change of Frequency (RoCoF) in case of a generation-demand imbalance \cite{Milano2018}. In the future, with more generation coming from wind and solar power, the ability of the system to maintain the frequency within the acceptable range is diminished. Indeed, photovoltaic systems are connected to the grid through inverters, which do not exhibit rotational inertia. Even in the case of inverter-interfaced wind generators, the inverter electrically decouples the rotor's rotational inertia from the system \cite{tielens2012grid}.

Transmission System Operators (TSOs) around the world are concerned with the stability issues associated with large penetration of renewable energy in their systems. In the United States, the Electric Reliability Council of Texas (ERCOT) has studied the effect of low inertia on the security and reliability of the grid \cite{ERCOTwind}. The Irish TSO, EirGrid, is designing ancillary services to remunerate providers of rotational or synthetic inertia \cite{EirGridReport}. Furthermore, EirGrid currently imposes limits on the maximum instantaneous penetration of variable RES with respect to the total load demand at any point in time.

In systems with low rotational inertia, TSOs must impose minimum inertia requirements in order to secure frequency stability and avoid system collapse in case of a severe fault or a sudden mismatch between generation and demand. With such new requirements, the traditional Unit Commitment (UC) problem, i.e., the day-ahead scheduling process to decide which generators will be committed, may be affected as more Synchronous Generators (SGs) could be dispatched for the sole purpose of providing inertia.

Several papers have approached the problem of including inertia requirements in the UC problem. In \cite{EnergiforskReport,perez2016robust,daly2015inertia}, the authors use the swing equation of Center-of-Inertia (CoI), which allows them to derive the RoCoF constraint and study its effect on the UC schedule. However, this approach oversimplifies the problem as it neglects metrics related to frequency deviation from the setpoint. This problem was addressed in \cite{Ahmadi} and \cite{Teng2016StochasticSW} with the inclusion of a constraint limiting the post-disturbance maximum frequency deviation (i.e., frequency nadir). In \cite{Ahmadi}, the analytic form of frequency nadir as a function of active power disturbance is derived using a system frequency model obtained from \cite{GenOrdModel}. The nadir expression is then linearized and added to the UC model, while considering a fixed sudden load increase. On the other hand, the authors of \cite{Teng2016StochasticSW} bypass the explicit modeling of turbine and governor control, as well as their impact on frequency dynamics, by imposing strict assumptions on system damping and total frequency response provision at each node. Moreover, they look at the impact of wind uncertainty on inertia requirements. While the simplifications proposed in \cite{Teng2016StochasticSW} enable the formulation of a nadir constraint without the explicit consideration of second-order frequency dynamics, they oversimplify the actual control implementation and disregard the aggregate impact of governor damping. Furthermore, none of the aforementioned studies incorporate the converter interface of RES and the impact of respective control schemes on the UC formulation.  

This study builds on the work of \cite{Ahmadi} and \cite{Teng2016StochasticSW}, and improves on it in several ways. First, we improve the frequency dynamics model in \cite{Ahmadi} by including the state-of-the-art converter control schemes of inverter-based generation, more specifically Virtual Synchronous Machine (VSM) and droop control. In contrast to the existing literature, where SG inertia and damping constants are usually numerically modified in order to compensate for high RES integration, we analyze a realistic model of a low-inertia system comprising both SG and converter dynamic models. This allows us to derive detailed analytic expressions of relevant frequency metrics as functions of multiple system variables (e.g., inertia, damping, aggregate droop gain, etc.) to be determined by the UC model, as opposed to the approach in \cite{Teng2016StochasticSW} where inertia constant was the only decision variable of interest. Moreover, in addition to frequency nadir and RoCoF, we incorporate the limitation on quasi steady-state frequency deviation into the UC formulation. Secondly, a more straightforward method is proposed to extract bounds for decision variables contributing to frequency nadir, which allows us to incorporate the non-linear nadir constraint in the UC problem in a more efficient way compared to \cite{Ahmadi}. Furthermore, similar to \cite{Teng2016StochasticSW}, this paper includes both the wind uncertainty and potential loss of generation in the UC model. However, we present a more comprehensive approach towards event probability computation and structuring of the scenario tree for the two-stage stochastic UC problem.

The rest of the paper is structured as follows. In Section~\ref{sec:FrequencyDynamics}, the derivation of post-contingency frequency dynamics in a low-inertia, multi-machine system is discussed. The obtained time-domain, analytic expressions are then linearized and incorporated into a stochastic UC formulation in Section~\ref{sec:FrequencyConstraints}. Subsequently, the modeling of uncertainties, namely equipment failure and wind power, in the form of probabilistic scenarios is presented in Section~\ref{sec:Uncertainty}. Section~\ref{sec:StochUC} provides the mathematical formulation of the stochastic UC problem. Finally, Section~\ref{sec:CaseStudy} presents and discusses the simulation results using a modified version of the IEEE RTS-96 system, whereas Section~\ref{sec:Conclusion} draws the main conclusions and discusses the outlook of the study.

\section{Low-Inertia System Frequency Dynamics} \label{sec:FrequencyDynamics}

\subsection{Inertial Response and Primary Frequency Control Model} \label{subsec:2.1}

\begin{figure}[b!] 
	\centering
	\vspace{-0.3cm}
	\scalebox{0.75}{\includegraphics[]{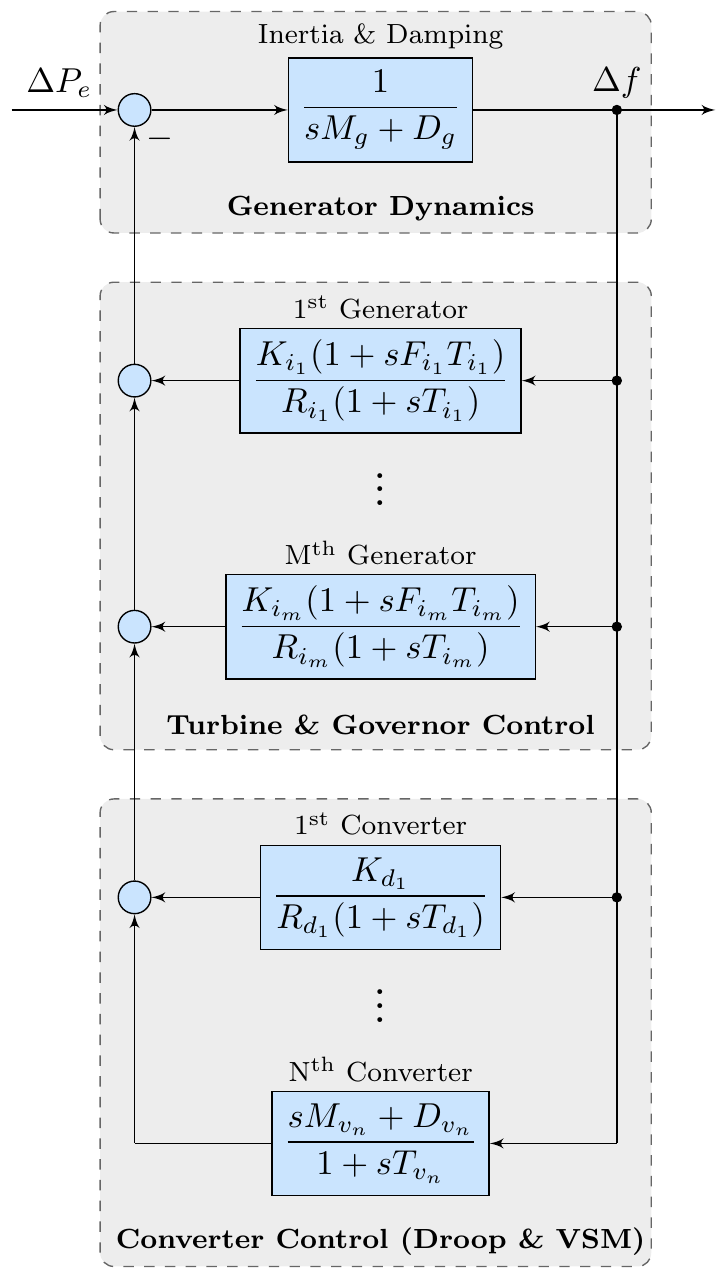}}
	\caption{Uniform system frequency dynamics model.}
	\label{fig:freq_dyn}
\end{figure}

We first focus on deriving a simplified, but sufficiently accurate, uniform frequency response model of a low-inertia system previously introduced in \cite{UrosLQR}. Let us consider a system comprised of \textit{traditional} $(i \in \mathcal{I})$ and \textit{converter-based} $(j \in \mathcal{J})$ generators depicted in Fig.~\ref{fig:freq_dyn}. 

The generator dynamics are described by the swing equation, with $M_g$ and $D_g$ denoting the normalized inertia and damping constants corresponding to the synchronous generators' CoI. The low-order model proposed in \cite{Anderson1990} is used for modeling the governor droop and turbine dynamics; $T_i$ are the turbine time constants, $R_i$ and $K_i$ are the respective droop and mechanical power gain factors, while $F_i$ refers to the fraction of total power generated by the turbines of synchronous machines. Furthermore, we incorporate the impact of grid-forming converters, as they are the only type of power electronic-interfaced units providing frequency support \cite{Rocabert2012,UrosGM}. A particular focus is set on droop $(d\in \mathcal{J}_d\subseteq\mathcal{J})$ and VSM $(v \in \mathcal{J}_v\subseteq\mathcal{J})$ control schemes, as two of the currently most prevalent emulation techniques in the literature, which in fact have equivalent properties in the grid-forming mode of operation \cite{Tamrakar2017}. Here, $T_d=T_v\equiv T_j$ are the time constants of all converters, $R_d$ and $K_d$ are the respective droop and electrical power gain factors, whereas $M_{v}$ and $D_{v}$ denote the normalized virtual inertia and damping constants of VSM converters.

\subsection{Analytic Derivation of Frequency Metrics} \label{subsec:2.2}
From Fig.~\ref{fig:freq_dyn} we can now derive a transfer function $G(s)$ of a general-order system dynamics, as follows:
\begin{align}
    G(s) &= \dfrac{\Delta f}{\Delta P_e} = \bigg(\underbrace{(sM_i+D_i )+\sum\limits_{i\in\mathcal{I}} \dfrac{K_i (1+sF_i T_i )}{R_i(1+sT_i)}}_{\text{traditional generators}} \nonumber \\
    & + \underbrace{\sum\limits_{d \in \mathcal{J}_d} \dfrac{K_d}{R_d (1+sT_d )}}_{\text{droop converters}} + \underbrace{\sum\limits_{v\in\mathcal{J}_v } \dfrac{sM_v+D_v} {1+sT_v} }_{\text{VSM converters}}\bigg)^{-1}. \label{eq:G1}
\end{align}
Assuming similar time constants $(T_i\approx T)$ of all synchronous machines, usually 2-3 orders of magnitude higher than the ones of converters, justifies the approximation $T \gg T_j\approx0$. Now we can transform \eqref{eq:G1} into the following expression:
\begin{equation}
    G(s) = \frac{1}{MT}\frac{1+sT}{s^2+2\zeta\omega_n s + \omega_n^2}, \label{eq:G2}
\end{equation}
where the natural frequency $(\omega_n)$ and damping ratio $(\zeta)$ are
\begin{equation}
    \omega_n = \sqrt{\frac{D+R_g}{MT}}, \quad \zeta = \frac{M+T(D+F_g)}{2\sqrt{MT(D+R_g)}}, \label{eq:wn}
\end{equation}
and parameters $(M,D)$ and $(F_g,R_g)$ represent weighted system and synchronous generator averages, respectively. More details on mathematical formulation can be found in \cite{UrosLQR}.


Assuming a stepwise disturbance in the electrical power $\Delta P_e(s) = -\Delta P/s$, we can derive the time-domain expression for frequency deviation as well as the nadir ($\dot{f}_\mathrm{max}$), RoCoF ($\dot{f}_\mathrm{max}$) and steady-state deviation ($\Delta f_\mathrm{ss}$) frequency metrics:
\begin{subequations} \label{eq:puConst}
\begin{align}
    \Delta f_\mathrm{max} &= - \frac{\Delta P}{D+R_g} \left( 1 + \sqrt{\dfrac{T(R_g-F_g)}{M}} e^{-\zeta\omega_n t_m} \right), \label{eq:nadir}\\
    \dot{f}_\mathrm{max} &= \dot{f}(t_0^+) = -\frac{\Delta P}{M}, \\
    \Delta f_\mathrm{ss} &= -\frac{\Delta P}{D + R_g},
\end{align}
\end{subequations}
with the introduction of the new variable $\omega_d = \omega_n\sqrt{1-\zeta^2}$. 

The accuracy of the proposed model has already been investigated and verified in \cite{UrosLQR}. We can conclude that the frequency metrics of interest are directly dependent on the average system parameters $M$, $D$, $R_g$ and $F_g$, and thus they could be regulated through the UC model. In particular, RoCoF and steady-state deviation can be explicitly controlled via $\dot{f}_{\mathrm{max}}\sim M^{-1}$ and $\Delta f_\mathrm{ss}\sim (D+R_g)^{-1}$, while nadir can be modeled using a highly non-linear function $\Delta f_{\mathrm{max}}\left(M,D,R_g,F_g\right)$.

\section{Formulation of Frequency Constraints} \label{sec:FrequencyConstraints}

The aforementioned frequency expressions in \eqref{eq:puConst} are incorporated as constraints into the stochastic UC problem, converted into SI and bounded by prescribed ENTSO-e thresholds \cite{entsoe}, as follows:
\begin{subequations}
\begin{align}
    &\left|\frac{f_b\Delta P}{D+R_g} \left( 1 + \sqrt{\dfrac{T(R_g-F_g)}{M}} e^{-\zeta\omega_n t_m} \right)\right| &\leq \Delta f_\mathrm{lim}, \label{nadir_constr} \\
    &\left|\frac{f_b\Delta P}{M}\right| \leq \dot{f}_\mathrm{lim}, \label{rocof_constr} \\
    &\left|\frac{f_b\Delta P}{D+R_g}\right| \leq \Delta f_\mathrm{ss,lim}, \label{qss_constr}
\end{align}
\end{subequations}
with $f_b=50\,\mathrm{Hz}$ being the base frequency; $\Delta f_\mathrm{lim}=0.4\,\mathrm{Hz}$ is the Under-Frequency Load Shedding (UFLS) trigger, while $\dot{f}_\mathrm{lim}=0.5\,\mathrm{Hz}/\mathrm{s}$ and $\Delta f_\mathrm{ss,lim}=0.2\,\mathrm{Hz}$ are the maximum permissible RoCoF and steady-state frequency deviation. 

Constraints \eqref{rocof_constr} and \eqref{qss_constr} are linear, unlike the non-linear frequency nadir constraint \eqref{nadir_constr}. In order to avoid the high computational burden of a Mixed-Integer Non-Linear Program (MINLP) formulation and have a measurable optimality gap, we use a linear approximation of \eqref{nadir_constr} which allows us to maintain a Mixed-Ineger Linear Program (MILP) formulation of the stochastic UC problem.

\subsection{Piece-wise Linearization of Nadir Expression} \label{subsec:PWL}

The study in \cite{Ahmadi} proposes a Piece-Wise Linearization (PWL) technique for obtaining a linearized expression for frequency nadir in order to subsequently integrate it into a UC problem. To improve clarity, this technique is outlined here and in Appendix~\ref{appendixA} before comparing its computational burden against our proposed approach introduced in Section~\ref{subsec:new_approach}. Let us recall from Section~\ref{subsec:2.2} that the frequency nadir expression is a function of four variables $(R_g,F_g,M,D)$, and as such too complicated to be directly handled by the PWL. Considering that the aggregate damping constant is of the form $D(D_i,D_v,R_d)$, with respective damping and droop gains usually strictly prescribed within narrow ranges by the system operator, it is justifiable to assume a constant $D$ and therefore becomes $\Delta f_\mathrm{max}(M,R_g,F_g)$. Hence, the PWL formulation aims to minimize the following objective function 
\vspace{-0.05cm}
\begin{align}
    \min_{\Psi} \sum_{\eta} \bigg(&\max_{1 \leq \nu \leq \overline{\nu} } \left\{a_\nu R_g^{(\eta)} + b_\nu F_g^{(\eta)} + c_\nu M^{(\eta)} + d_\nu \right\} \nonumber \\ 
    & - \Delta f_\mathrm{max}\left(R_g^{(\eta)}, F_g^{(\eta)}, M^{(\eta)}\right)\bigg)^2 \label{eq:3d_pwl_obj}
\end{align}
with $\Psi=\{a_\nu,b_\nu,c_\nu,d_\nu, \forall \nu\}$, being the set of optimization variables, $\eta$ denoting the evaluation point and $\nu$ referring to the number of PWL segments. The objective function \eqref{eq:3d_pwl_obj} penalizes the difference between the appropriate PWL segment and the nadir function at all evaluation points. Given the convex nature of the nadir function, the inner $\max$ operator chooses the appropriate PWL segment for each evaluation point by looking at which segment is closest to the curve at that specific point. To improve understanding, we provide an illustration of the PWL method from \cite{Ahmadi} on a simple one variable function in Appendix~\ref{appendixA}.

\begin{figure}[!h] 
	\centering
	\vspace{-0.3cm}
	\scalebox{0.425}{\includegraphics[]{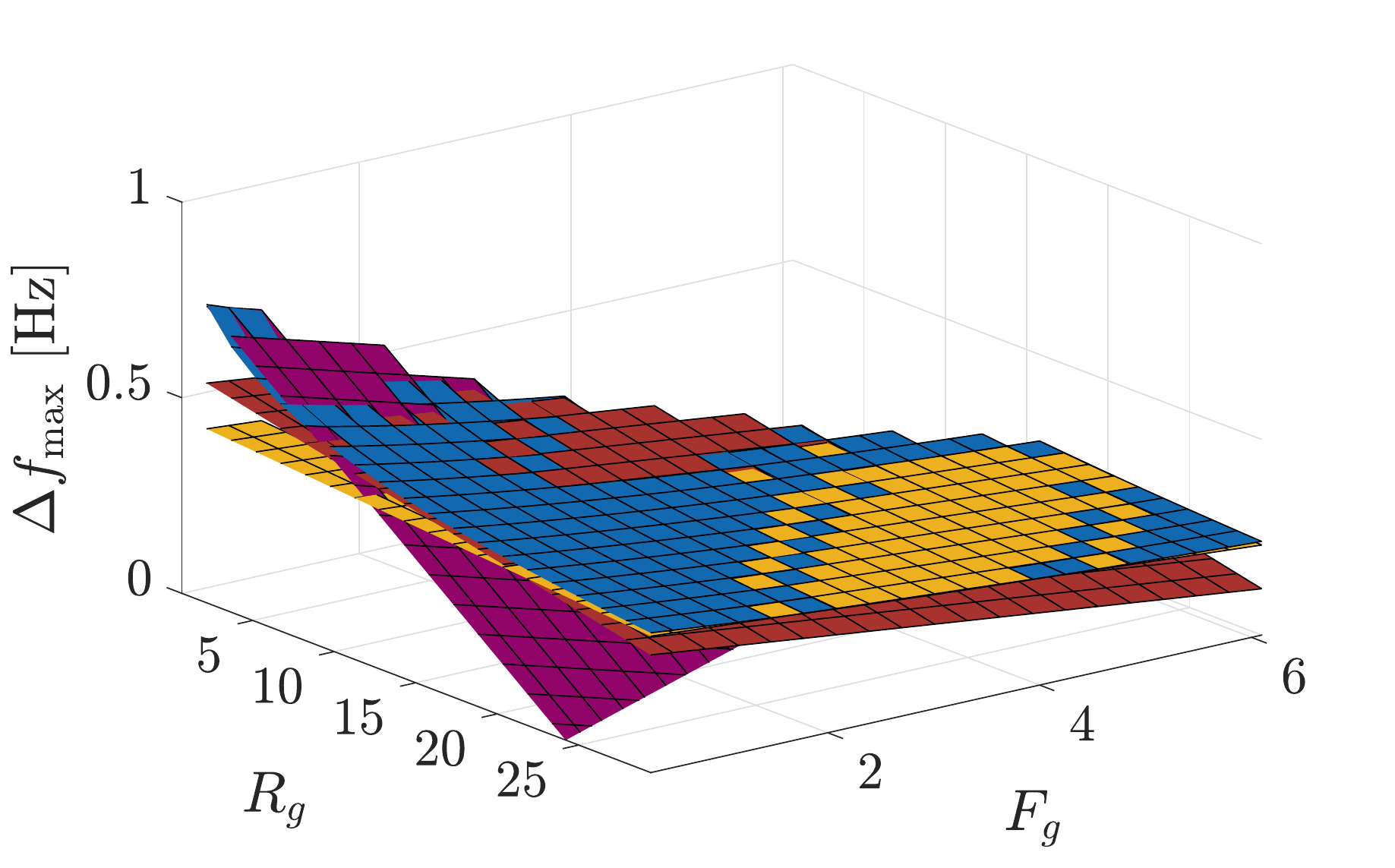}}
	\caption{PWL of the nadir constraint for $M=9$.}
	\label{fig:pwl3d}
\end{figure}

Upon obtaining the optimal solution of the model in \eqref{eq:3d_pwl_obj}, denoted as $(a^*_\nu,b^*_\nu,c^*_\nu,d^*_\nu)$, the nadir constraint can be integrated into the MILP UC model by adding a set of inequalities described in Appendix~\ref{appendixA}, along with the nadir threshold constraint of the form $f_b\,t_3 \leq \Delta f_\mathrm{lim}$.
The results for the approximation of frequency nadir function for a test system of 20 generators described in Section~\ref{sec:CaseStudy} are shown in Fig.~\ref{fig:pwl3d}, where a loss of the largest unit is considered. Note that Fig.~\ref{fig:pwl3d} showcases the surface plot for a fixed inertia constant $M$ and thus ignores one degree of freedom. The original surface is presented in blue, whereas its PWL-approximation segments are the planes depicted in various colors. It is important to note that the optimization problem \eqref{eq:3d_pwl_obj} is computationally intensive and thus in order to obtain results within reasonable computational time, the number of PWL segments used for the approximation as well as the number of evaluation points have to be kept low.


\subsection{Extracting Bounds on Relevant Variables} \label{subsec:new_approach}

\begin{figure}[!b]
    \centering
    \vspace{-0.3cm}
    \scalebox{0.425}{
    \includegraphics[trim={1.8cm, 9cm, 2cm, 9.5cm},clip=true]{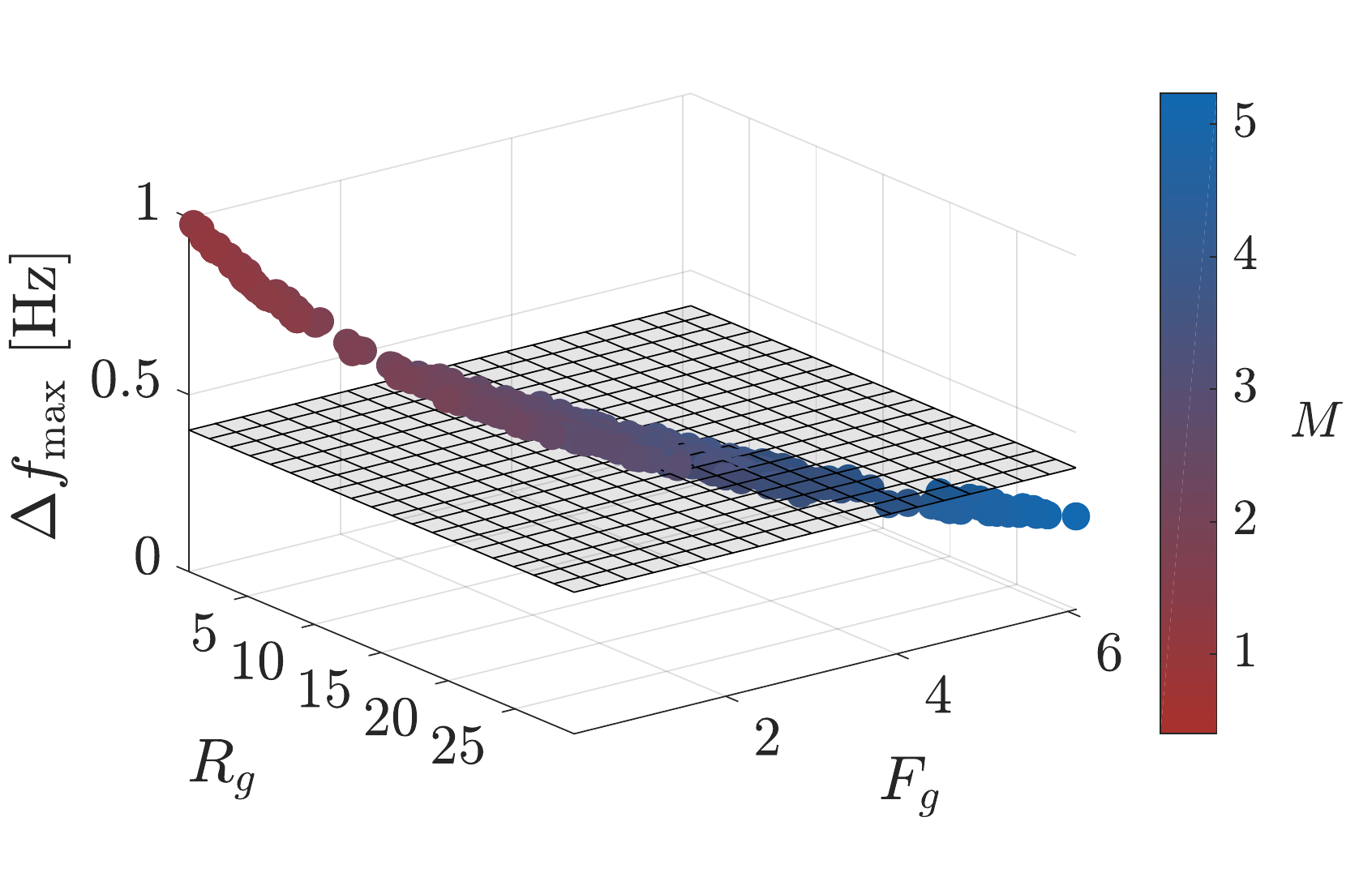}
    }
    \caption{All possible values of the nadir after a generator loss.}
    \label{fig:nadir_3d_scatter}
\end{figure}

An alternative approach for linearizing the nadir constraint and integrating it into the UC problem is to confine the values of $R_g$, $F_g$, $M$ and $D$ within a plausible range to guarantee that the nadir threshold in \eqref{nadir_constr} is not violated. With this approach the damping variable $D$ can easily be included and does not need to be set constant.
The scatter plot presented in Fig.~\ref{fig:nadir_3d_scatter} reflects all possible values of frequency nadir after the loss of the largest generator, for the same system as in Section~\ref{subsec:PWL}. In the general case, for a system that comprises $|\mathcal{I}|$ generators, there will be $2^{|\mathcal{I}|-1}$ possible generator commitment combinations after a generator outage. By obtaining the set of these dispatch combinations, the values of $R_g$, $F_g$ and $M$ at which the UFLS threshold is not violated can be extracted, corresponding to the points below the shaded plane in Fig.~\ref{fig:nadir_3d_scatter}. Subsequently, these values are used to substitute the nadir constraint in the unit commitment as follows:
\begin{equation}
    F_g \geq F_g^{\mathrm{lim}},\quad R_g \geq R_g^{\mathrm{lim}},\quad M \geq M^{\mathrm{lim}}.
    \label{eq:limits_nadir}
\end{equation}

Table~\ref{tab:comp_times} provides a comparison of the proposed method to the PWL technique, in terms of computational time that is needed to obtain the equivalent linear nadir equations for a single value of $\Delta P$. It is clear that the PWL is more time intensive, especially when aiming for an increased precision. It should be noted though that for very large systems, the calculations of $2^{|\mathcal{I}|-1}$ combinations for the bound extraction method would become more computationally expensive. Indeed, the computation time increases by a factor of $2^{\Delta |\mathcal{I}|}$ for every additional $\Delta |\mathcal{I}|$ generators included in the system. For the purposes of this paper, the proposed bound extraction method will be used as it is significantly faster and introduces less error when applied to the 20-generator test system under investigation.

\begin{table}[!t]
\renewcommand{\arraystretch}{1.2}
\caption{Computational cost of the linearization methods.}
\label{tab:comp_times}
\noindent
\centering
    \begin{minipage}{\linewidth} 
    \renewcommand\footnoterule{\vspace*{-5pt}} 
    \begin{center}
        \begin{tabular}{ c | c }
            \toprule
            \textbf{Linearization technique} & \textbf{Computational time~[s]} \\
            \cline{1-2}
            PWL ($m=3$,\,$k=4$) &  70\\
            \arrayrulecolor{black!30}\hline
            PWL ($m=4$,\,$k=4$) &  7200\\
            \arrayrulecolor{black!30}\hline
            Bound extraction & 20 \\
            \arrayrulecolor{black}\bottomrule
        \end{tabular}
        \end{center}
    \end{minipage}
    \vspace{-0.3cm}
\end{table}

\section{Modeling equipment-failure and wind power uncertainties}  \label{sec:Uncertainty}

This section describes the modeling of uncertainty pertaining to equipment failure and wind power production during power system operation. The uncertain nature of wind power production is modeled using a set of scenarios $\Omega$ that captures the spatio-temporal interdependence of forecast errors, for every wind farm location and during the whole scheduling horizon. Each wind power realization scenario $\omega$ has the same probability of occurrence denoted as $\pi_{\omega}$. 

In terms of equipment failure uncertainty, we consider as the set of credible contingencies, the unforeseen outages of synchronous generators, whereas transmission assets are assumed to be $100\,\%$ reliable. 
In order to reduce the computational burden, we follow the assumption from \cite{AhmadiKhatir2013MultiAreaEA} considering that the generation outages happen at a discrete time period, while failed assets remain unavailable for the rest of the scheduling horizon, i.e., the Mean Time To Repair (MTTR) is greater than the scheduling horizon of the day-ahead electricity market.  
For the purpose of assessing the impact of frequency constraints on the unit commitment schedule, we consider as contingency period the one in which the power system faces the highest wind power penetration as scarcity of inertia is most likely to occur in this time due to the displacement of synchronous generators from the day-ahead schedule.

To calculate the probability $\pi_c$ associated with contingency scenario $c$, we index the set of credible contingencies by $\kappa=1,2,\dots,\mathcal{K}$. We denote by $A(\kappa,\tau)$ the random event of contingency $\kappa$ happening within time period $\tau$. Random event $B(\kappa)$ corresponds to contingency $\kappa$ not occurring during the entire scheduling horizon.
We further denote as $\lambda_\kappa$ the inverse of the Mean Time To Failure $(\mathrm{MTTF})$ $\kappa$, i.e., $\lambda_\kappa = 1/\mathrm{MTTF}$.
Considering we are looking at only one hour in which the outages may occur, the probability $\pi_c$ for each contingency scenario $\kappa$ is derived from the probabilities of occurrence of random events $A(\kappa,\tau)$ and $B(\kappa)$ that are calculated using the following expressions according to \cite{conejo2010decision}:
\begin{subequations}
\label{eq:ContProb}
\begin{align}
    & \pi[A(1,\tau)] = \mathrm{exp}(-\lambda_1 \tau)(\mathrm{exp}(\lambda_1)\!-\!1),  \\
    & \pi[A(1,\tau)] =  \pi[A(2,\tau)] = \ldots = \pi[A(\mathcal{K},\tau)], \\
    & \pi[B(1)] = \mathrm{exp}(-\lambda_1 \tau),  \\
    & \pi[B(1)] = \pi[B(2)] = \ldots = \pi[B(\mathcal{K})].
\end{align}
\end{subequations}
Assuming statistical independence between all contingencies, the probability $\pi_{c_0}$ of the no-contingency scenario is equal to
\begin{align}
    \pi_{c_0} &= \prod_{\kappa=1}^{\mathcal{K}} \pi[B(\kappa)], \label{eq:ProbNC}
\end{align}
while the probability $\pi_{c_\kappa}$ of losing a generator is equal to
\begin{equation}
    \begin{aligned}
    \pi_{c_\kappa} = \pi[A(1,\tau)] \prod_{\substack{y=1 \\ y \neq k }}^{\mathcal{K}} \pi[B(y)], \quad \forall \kappa = 1,...,\mathcal{K}. \label{eq:ProbCk}
    \end{aligned}
\end{equation}
It should be noted that the sum of probabilities $\pi_{c_0}$ and $\pi_{c_\kappa}$ is lower than 1, since sequential contingencies are not considered. For instance, setting $\mathrm{MTTF}$ equal to $1000\,\mathrm{h}$ for all generators, we obtain from \eqref{eq:ContProb}  $\pi[A(\kappa,\tau)] = 0.9995\!\times\!10^{-3}$ and  $\pi[B(\kappa)] = 0.9990$. According to \eqref{eq:ProbNC} and \eqref{eq:ProbCk} we obtain respectively $\pi_{c_0} = 0.9960 $ and $\pi_{c_\kappa} = 0.9965 \times 10^{-3}$ and thus $\sum_{\kappa=0}^{4} \pi_{c_\kappa} = 0.9999 \approx 1$.

\begin{figure}[!b]
    \centering
    \vspace{-0.3cm}
    \scalebox{0.69}{
    \includegraphics[]{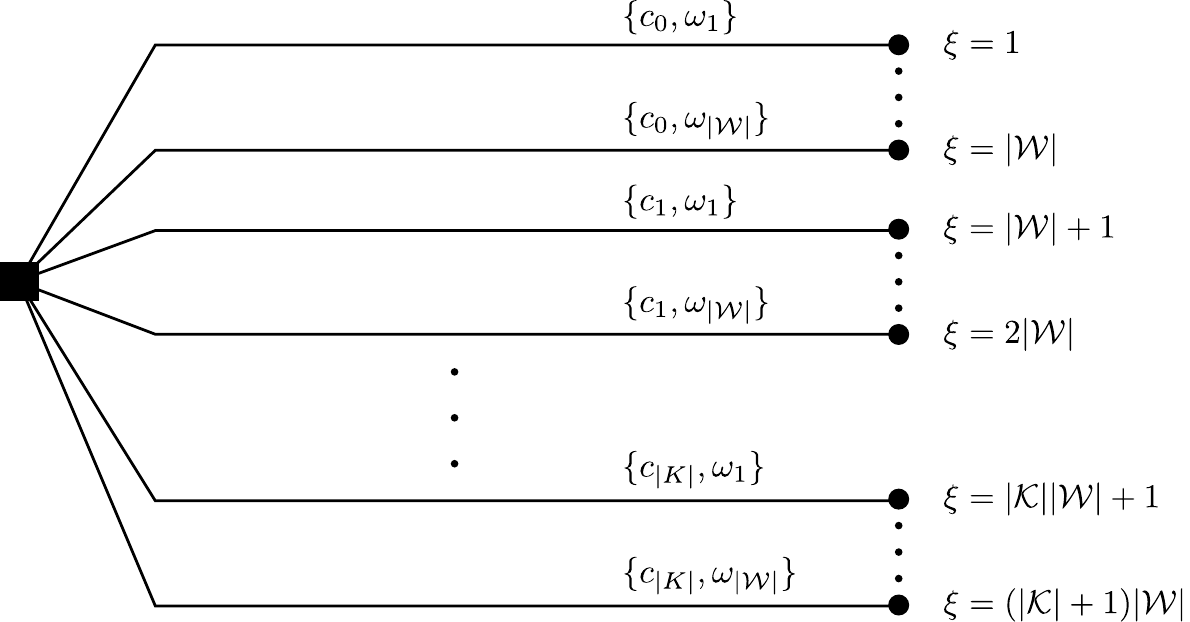}
    }
    \caption{Scenario tree for the two-stage stochastic UC problem.}
    \label{fig:scenario_tree}
\end{figure}

Combining the scenarios modeling the equipment failure and wind power uncertainty into a single scenario set $\mathcal{E}$, we define each scenario $\xi$ as a pair of contingency $c$ and wind power realization $\omega$.
For each $\xi = \{c,\omega\}$ the corresponding probability of occurrence is given as $\pi_{\xi} = \pi_{\omega} \cdot \pi_{c}$ and $\sum_{\xi \in \mathcal{E}} \pi_{\xi} \approx 1$, assuming that equipment outages and wind power production are statistically independent events. The structure of the scenario set $\mathcal{E}$ used in the stochastic UC formulation is illustrated as the scenario tree shown in Fig.~\ref{fig:scenario_tree}, for $|\mathcal{K}|$ contingencies and $|\mathcal{W}|$ wind power scenarios.

\section{Stochastic Unit Commitment} \label{sec:StochUC}

This section provides the mathematical formulation of the stochastic unit commitment \cite{delikaraoglou2014high}, with an addition of frequency-related constraints. The proposed model is a two-stage stochastic optimization problem which can be written as:
\vspace{-0.15cm}
\begin{subequations}
\begin{equation}
    \begin{aligned}
        &\min_\Phi \,\, \sum_{t\in \mathcal{T}} \sum_{i\in \mathcal{I}} \big(C^{\mathrm{SU}}_{i}y_{it} + C^{\mathrm{SD}}_{i}z_{it} + C_{i} p_{it}\big)\ + \\
        & \sum_{t\in \mathcal{T}} \sum_{\xi \in \mathcal{E}} \pi_{\xi} \Big[\sum_{i\in \mathcal{I}} \big(C^{+}_{i}r^{+}_{i\xi t} - C^{-}_{i}r^{-}_{i\xi t} \big) + \sum_{n\in \mathcal{N}}C^{\mathrm{sh}}\,l^{\mathrm{shed}}_{n\xi t}\Big] \label{eq:UC_obj}
    \end{aligned}
\end{equation}
subject to
\begin{flalign}
    &\sum_{i \in \mathcal{I}_n} p_{it}\! + \! \sum_{j \in \mathcal{J}_n } w_{jt} - D_{nt} -  \nonumber \\ 
    & \quad \sum_{m: (n,m) \in \mathcal{L}} B_{nm} (\hat{\delta}_{n t} - \hat{\delta}_{m t}) = 0 , \quad \forall n, t, \label{eq:pow_bal_da}\\
	& B_{nm} (\hat{\delta}_{nt} - \hat{\delta}_{mt}) \le \overline{f}_{nm}, \quad \forall (n,m) \in \mathcal{L},t, \label{eq:flowlim_da}\\
    &u_{it} - u_{i(t-1)} \leq u_{i\tau_{i}^1}, \quad \forall i,t,  \label{eq:min_on}\\ 
    &u_{i(t-1)} - u_{it} \leq 1 - u_{i\tau_{i}^0}, \quad \forall i,t, \label{eq:min_off}\\
	&y_{it} \geq u_{it} - u_{i(t-1)},\quad \forall i,t, \label{eq:startup}\\
	&z_{it} \geq u_{i(t-1)} - u_{it},\quad \forall i,t, \label{eq:shutdown}\\
	& \sum_{i \in \mathcal{I}_n} \Big[ r^{+}_{i\xi t} - r^{-}_{i\xi t} - p_{it}(1-\alpha_{i\xi t}) \Big]  \nonumber\\ 
	& + \sum_{m: (n,m) \in \mathcal{L}} B_{nm} (\hat{\delta}_{n t} - \tilde{\delta}_{n \xi t} - \hat{\delta}_{m t} + \tilde{\delta}_{m \xi t})  \label{eq:pow_bal_rt}\\ 
	& + \sum_{j \in \mathcal{J}_n} (W^{*}_{j\xi t} - w_{jt} - w^{\mathrm{spill}}_{j\xi t}) + l^{\mathrm{shed}}_{n \xi t}= 0,  \quad \forall n, \xi, t, \nonumber\\
	&p_{it} + r^{+}_{i\xi t} \leq \overline{P_{i}} u_{it}, \quad \forall i,\xi, t, \label{eq:pmax}\\
	&p_{it} - r^{-}_{i\xi t} \geq \underline{P_{i}} u_{it}, \quad \forall i,\xi, t,  \label{eq:pmin}\\
	&p_{it} - p_{i(t-1)} + r^{+}_{i\xi t} - r^{+}_{i\xi (t-1)} \leq R^{\mathrm{U}}_{i}, \quad \forall i,\xi,t, \label{eq:rampup}\\
	&p_{it}\! -\! p_{i(t-1)} - r^{-}_{i\xi t}\! +\! r^{-}_{i\xi (t-1)} \geq -R^{\mathrm{D}}_{i},  \quad \forall i,\xi,t, \label{eq:rampdown} \\
	&r^{+}_{i\xi t} \leq R_i^{+} \alpha_{i\xi t}, \quad \forall i,\xi, t,  \label{eq:rmax}\\
	&r^{-}_{i\xi t} \leq R_i^{-} \alpha_{i\xi t}, \quad \forall i,\xi, t,  \label{eq:rmin}\\
	& B_{nm} (\tilde{\delta}_{nt} - \tilde{\delta}_{mt}) \le \overline{f}_{nm}, \quad \forall (n,m) \in \Lambda,t, \label{eq:flowlim_rt}\\
	& w^{\mathrm{spill}}_{j\xi t} \leq w_{j\xi t}, \quad \forall j,\xi,t,  \label{eq:spill}\\
	& l^{\mathrm{shed}}_{n\xi t} \leq D_{n t}, \quad \forall n,\xi,t,  \label{eq:Lshed}\\
	& k_{i\xi t} = \frac{P_{i}  K_i }{\sum_{i\in \mathcal{I}} P_{i} + W_{v}+W_{d}}  u_{it} \alpha_{i\xi t}, \quad \forall i,\xi,t, \label{eq:iner_K}\\
	&F_{\xi t} = \sum_{i\in \mathcal{I}} \frac{F_i k_{i\xi t}}{R_i}, \quad \forall t,\,\xi, \label{eq:iner_F}\\
	&R_{\xi t} = \sum_{i\in \mathcal{I}} \frac{k_{i\xi t}}{R_i},  \quad \forall t,\,\xi, \label{eq:iner_R}\\
	&M_{\xi t} = \sum_{i\in \mathcal{I}} 2 H_g k_{i\xi t},  \quad \forall t,\,\xi, \label{eq:iner_M}\\
	&\frac{\dot{f}_\mathrm{lim}}{f_b} (M_{\xi t}+M_v) \geq \Delta P_{\xi t},  \quad \forall t,\,\xi, \label{eq:rocof_cstr} \\
	&F_{\xi t} \geq F_{\xi t}^{\mathrm{lim}},\ R_{\xi t} \geq R_{\xi t}^{\mathrm{lim}},\ M_{\xi t} +M_v \geq M_{\xi t}^{\mathrm{lim}},  \quad \forall t,\,\xi, \label{eq:nadir_cstr}\\
    &\frac{\Delta f_\mathrm{ss,lim}}{f_b} (D+R_{\xi t}) \geq \Delta P_{\xi t},  \quad \forall t,\,\xi, \label{eq:qss_cstr} \\
    & p_{it} \ge 0, \forall i,t ;\, w_{jt} \ge 0, \forall j,t; \, \hat{\delta}_{n t} \ge 0, \forall n,t; \,k_{i\xi t} \ge 0, \forall i ,\xi , t; \nonumber \\
    & r^{+}_{i \xi t}, r^{-}_{i \xi t} \ge 0, \forall i,\xi, t ; \, l^{\mathrm{shed}}_{n\xi t} \ge 0, \forall n, \xi, t ; \, w^{\mathrm{spill}}_{j \xi t} \ge 0, \; \forall j, \xi, t; \nonumber \\
    & F_{\xi t},R_{\xi t},M_{\xi t} \ge 0, \forall \xi, t;\, u_{it},y_{it},z_{it} \in \{0,1\}, \label{eq:vardef} 
\end{flalign}
\end{subequations}
where $\Phi = \{ p_{it},\, u_{it},\, y_{it},\, z_{it},\, \forall i, t; \; w_{jt}, \forall j, t; \; \hat{\delta}_{nt}, \forall n, t; $
$\tilde{\delta}_{n \xi t}, \forall n, \xi, t ; \, r^{+}_{i\xi t},\, r^{-}_{i\xi t}, \forall i, \xi, t ; \,  w^{\mathrm{spill}}_{j \xi t}, \forall j, \xi , t ; \, l^{\mathrm{shed}}_{n \xi t}, \forall n ,\xi , t; \,$ $k_{i\xi t}, \forall i ,\xi , t ; \,  F_{\xi t},\, R_{\xi t},\, M_{\xi t}, \forall \xi, t \}$ is the set of optimization variables.

The objective function \eqref{eq:UC_obj} to be minimized is the total expected system cost that comprises the day-ahead energy and the real-time balancing costs. The day-ahead component consists of the fuel costs $C_{i}$ as well as the start-up and shut-down costs. 
The real-time component includes the re-dispatch cost from the deployment of upward and downward reserves based on the corresponding offer prices $C^{+}_{i}$ and $C^{-}_{i}$, as well as the involuntary load shedding at the value of lost load $C^{\mathrm{sh}}$.

Equation \eqref{eq:pow_bal_da} enforces the nodal power balance of the day-ahead schedule, while network power flows at the day-ahead stage are restricted by the transmission capacity limits in \eqref{eq:flowlim_da}.
Constraints \eqref{eq:min_on}-\eqref{eq:min_off} model the minimum online and offline time of conventional units based on commitment variable $u_{it}$, where parameters $\tau^1_i$ and $\tau^0_i$ are defined as $\tau_i^1 = \min \{t+T^1_i-1,\ T\}$ and $\tau_{i}^0 = \min \{t+T^0_i-1,\ T\} $, and $T^1_i$ and $T^0_i$ denote the duration that unit $i$ should remain online and offline, respectively. Constraints \eqref{eq:startup}-\eqref{eq:shutdown} model the start-up and shut-down of conventional units using the binary variables $y_{it}$ and $z_{it}$, respectively.
The real-time power balance for every uncertainty realization $\xi$ is enforced by constraint \eqref{eq:pow_bal_rt}. Parameter $\alpha_{i\xi t}$ models the availability of the generators to provide reserves, i.e., $\alpha_{i\xi t}$ is equal to 1 if generator $i$ at scenario $\xi$ and time $t$ is online and able to provide reserves and zero otherwise.
The scheduled energy production and the deployment of upward ($r^{+}_{i\xi t}$) and downward ($r^{-}_{i\xi t}$) reserves in each scenario $\xi$ are bounded by the generation capacity limits of each unit by constraints \eqref{eq:pmax}-\eqref{eq:pmin}, whereas constraints \eqref{eq:rampup}-\eqref{eq:rampdown} enforce the upward and downward ramping limits accounting for the real-time reserve activation.
Constraints \eqref{eq:rmax}-\eqref{eq:rmin} account for the limits of reserve capacity offers. 
Transmission capacity limits during real-time operation are enforced by constraint \eqref{eq:flowlim_rt}, whereas wind spillage $w^{\mathrm{spill}}_{j \xi t}$ and load shedding $l^{\mathrm{shed}}_{n \xi t}$ are bounded by the wind power realization and the nodal demand through constraints \eqref{eq:spill} and \eqref{eq:Lshed}, respectively. 

The set of constraints \eqref{eq:iner_K}-\eqref{eq:nadir_cstr} models the frequency limits of the power system. The equality constraint \eqref{eq:iner_K} defines $k_{i\xi t}$ as the gain factor $K_i$ of generator $i$ scaled by the ratio of its capacity over the total system capacity, which in turn is multiplied by the binary variable $u_{it}$ and the parameter $\alpha_{i\xi t}$ to indicate that a unit can only provide inertial response if it is committed and does not face an outage. Similarly, constraints \eqref{eq:iner_F}-\eqref{eq:iner_M} define average system variables for power fraction, droop and inertia, respectively. Constraint \eqref{eq:rocof_cstr} enforces the RoCoF limit, while nadir equivalent and quasi steady-state frequency bounds are imposed by constraints \eqref{eq:nadir_cstr} and \eqref{eq:qss_cstr}. Finally, constraints \eqref{eq:vardef}  are variable declarations.

\section{Case Study} \label{sec:CaseStudy}

\subsection{System Description}

\begin{figure}[!t] 
	\centering
	\scalebox{0.185}{\includegraphics[]{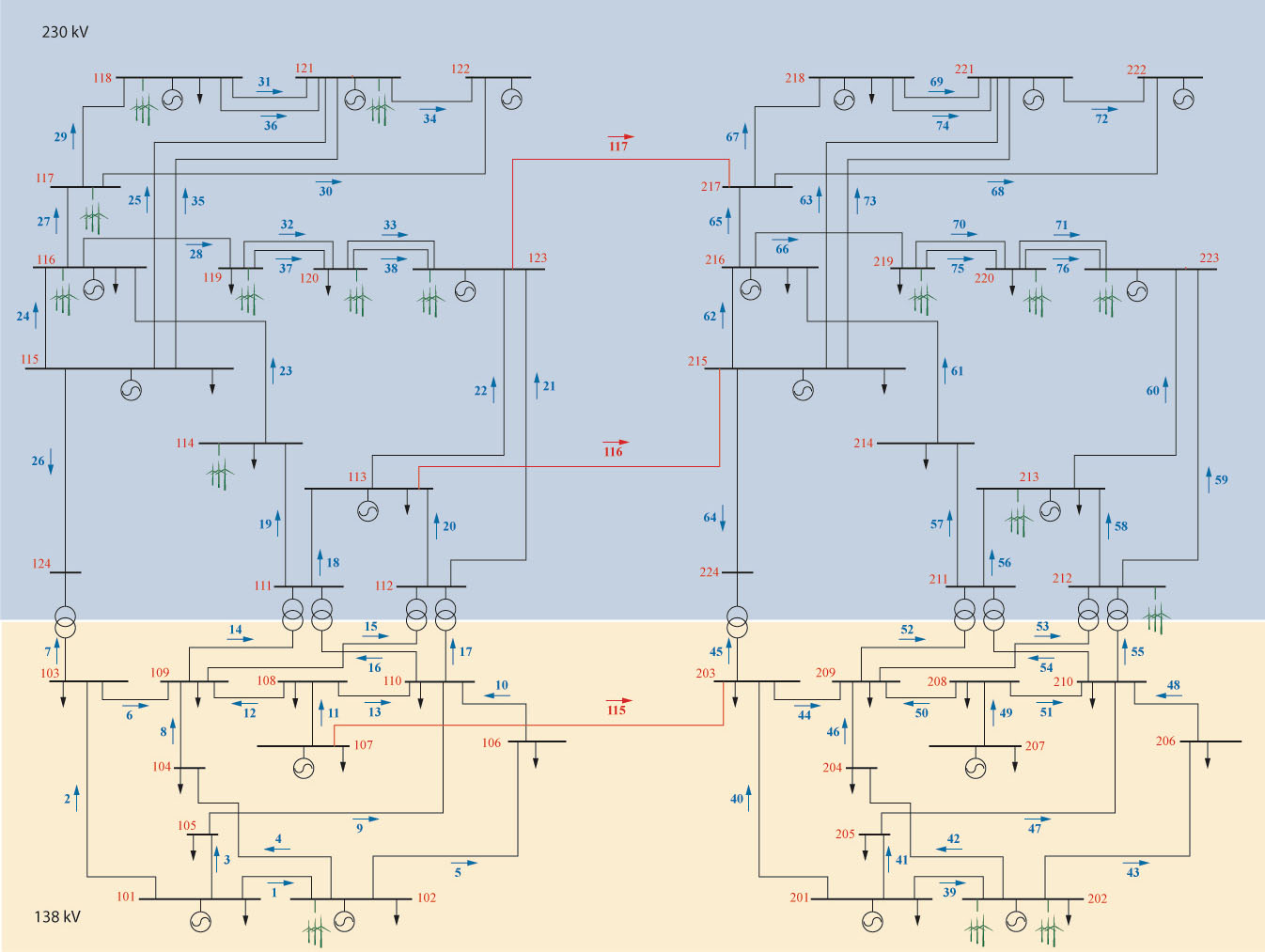}}
	\caption{Modified IEEE RTS-96 system diagram comprised of 2 areas, 16 wind farms and 20 synchronous generators \cite{modifiedRTS96}.}
	\label{fig:RTS_96}
	\vspace{-0.3cm}
\end{figure}

In order to analyze the performance of the stochastic UC model presented in Section~\ref{sec:StochUC}, we investigate a modified version of the IEEE RTS-96 power system from \cite{modifiedRTS96} depicted in Fig.~\ref{fig:RTS_96}, with 48 buses comprising areas 1 and 2 of the original system. Table \ref{tab:gen_param} shows the relevant parameters of different thermal plant types. The studied system includes 20 generators and 16 wind farms. We assume that six wind farms are providing virtual inertia; four of them via VSM control, and the remaining two through equivalent droop regulation. The UC is ran for two days without frequency constraints in order to initialize the system prior to introducing the frequency constraints on days 3, 4 and 5. This is done to ensure the impact of start-up costs are well distributed and not concentrated on one day. Therefore, the total simulation horizon is five days ($T=120\,\mathrm{h}$) whereas the UC schedule is optimized separately for each day, with the last hour of each day used as an input for the next.

As the set of possible contingencies we consider the failure of synchronous generators $i=\{1,6,8,10\}$. These generators are of various capacities, ranging from the smallest to the largest unit in the system. The hour 19 of day 3 (i.e., $t=67\,\mathrm{h}$) is selected to be the time instance of a possible contingency, as this is the hour with high wind penetration and low demand. Ten wind power scenarios are considered, which brings the total number of scenarios to 50. The optimization problem is formulated in Python and uses the Gurobi solver with default parameterization.

\begin{table}[!h]
\renewcommand{\arraystretch}{1.2}
\caption{Parameters of the thermal plants and VSM.}
\label{tab:gen_param}
\noindent
\centering
    \begin{minipage}{\linewidth} 
    \renewcommand\footnoterule{\vspace*{-5pt}} 
    \begin{center}
        \begin{tabular}{ c || c | c | c | c | c }
            \toprule
            \textbf{Type} & $H_g\,[\mathrm{s}]$ & $K_g\,[\mathrm{p.u.}]$ & $F_g\,[\mathrm{p.u.}]$ & $R_g\,[\mathrm{p.u.}]$ & $D_g\,[\mathrm{p.u.}]$ \\ 
            \cline{1-6}
            Nuclear & 4.5 & 0.98 & 0.25 & 0.04 & 0.6 \\
            \arrayrulecolor{black!30}\hline
            CCGT & 7.0 & 1.1  & 0.15 & 0.01 & 0.6 \\
            \arrayrulecolor{black!30}\hline
            OCGT & 5.5 & 0.95 & 0.35 & 0.03 & 0.6 \\
            
            \arrayrulecolor{black!30}\hline
            VSM & 6.0 & 1.0 & - & - & 0.6 \\
            \arrayrulecolor{black!30}\hline
            Droop & - & 1.0 & - & 0.05 & - \\
            \arrayrulecolor{black}\bottomrule
        \end{tabular}
        \end{center}
    \end{minipage}
\end{table}

\subsection{Results}

\begin{figure}[!t] 
	\centering
	\scalebox{1.1}{\includegraphics[]{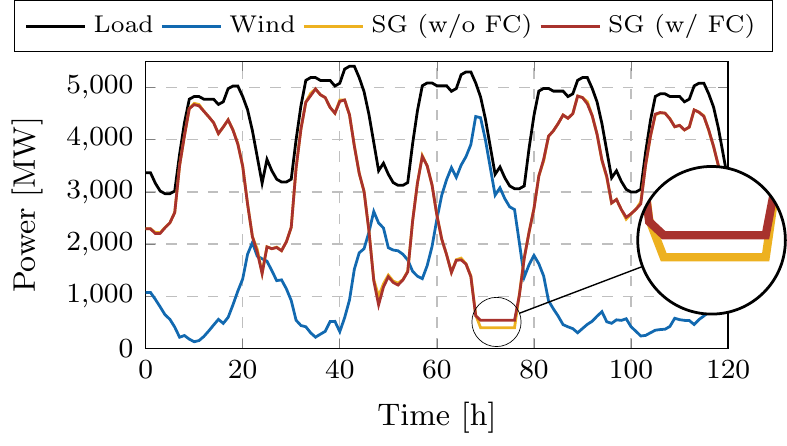}}
	\caption{UC dispatch of synchronous generation for respective load and wind profiles.}
	\label{fig:UC_profile}
	\vspace{0.05cm}
\end{figure}

In this section, the simulation results from the stochastic UC are presented. Fig.~\ref{fig:UC_profile} showcases the load and wind power profiles as well as the aggregate dispatch of synchronous generation for two UC runs: (i) without frequency constraints; and (ii) with frequency constraints. Furthermore, Table \ref{tab:commitment_comp} indicates the difference in the total number of generators committed between the two runs. Both Fig.~\ref{fig:UC_profile} and Table \ref{tab:commitment_comp} suggest that, although the amount of committed generators increases significantly, the total SG production is only slightly changed. This is justified by the fact that the extra generators are solely committed for the purpose of providing inertia, and are thus operating at their technical minimum. The production surplus arising from the additionally committed units is compensated by wind curtailment and other generators reducing their power output.

\begin{table}[!t]
\renewcommand{\arraystretch}{1.2}
\caption{Comparison of the total number of dispatched generators through UC for each hour.}
\label{tab:commitment_comp}
\noindent
\centering
    \begin{minipage}{\linewidth} 
    \renewcommand\footnoterule{\vspace*{-5pt}} 
    \begin{center}
        \begin{tabular}{ c || c | c | c | c | c | c | c | c | c }
            \toprule
            \textbf{Hour} & \textbf{65} & \textbf{66} & \textbf{67} & \textbf{68} & \textbf{69} & \textbf{70} & \textbf{71} & \textbf{72} & \textbf{73} \\
            \cline{1-10}
            w/o FC & 6  & 5  & 4  & 4  & 4  & 4  & 4  & 4  & 4   \\
            \arrayrulecolor{black!30}\hline
            w/ FC & 6 & 5 & \cellcolor{myGreen!20}10 & \cellcolor{myGreen!20}10 &	\cellcolor{myGreen!20}10 & \cellcolor{myGreen!20}10 & \cellcolor{myGreen!20}10 & \cellcolor{myGreen!20}10 & \cellcolor{myGreen!20}10 \\
            \arrayrulecolor{black}\bottomrule
        \end{tabular}
        \end{center}
    \end{minipage}
    \vspace{-0.3cm}
\end{table}

\begin{figure}[!b] 
	\centering
	\vspace{-0.3cm}
	\scalebox{1.1}{\includegraphics[]{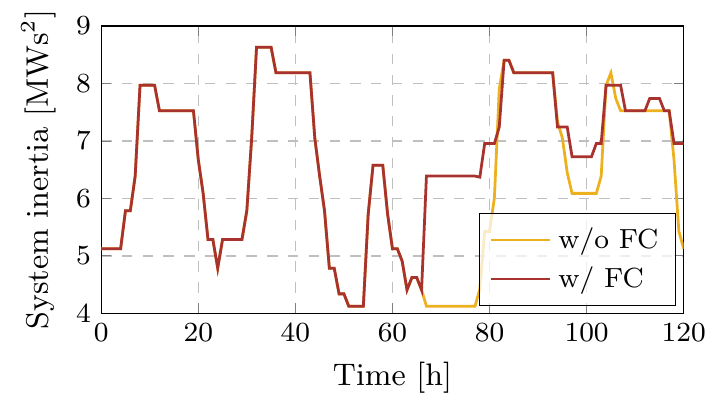}}
	\caption{Impact of frequency constraints on the aggregate level of system inertia.}
	\label{fig:system_inertia}
\end{figure}

Moreover, the evolution of aggregate system inertia over the course of the whole scheduling horizon is depicted in Fig.~\ref{fig:system_inertia}. A noticeable step change in total system inertia at hour $67$ reflects the violation of frequency constraints under contingency, which subsequently triggers a dispatch of auxiliary synchronous generators. While the inertia levels do not differentiate between the two scenarios during the first two days, on the days following the potential outage some carryover impacts can be observed. This is a consequence of the commitment schedule being radically changed at hour $t=67$, thus affecting the UC schedules in the following days.

\begin{figure}[!t]
    \centering
    \scalebox{1.1}{\includegraphics[]{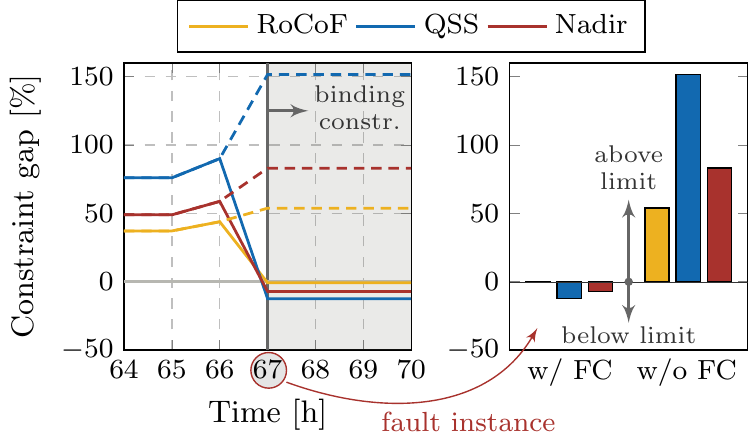}}
    \caption{Constraint gaps for different frequency metrics. Dashed lines refer to the scenario without frequency constraints.}
    \label{fig:IC_lim_values}
    \vspace{-0.3cm}
\end{figure}

\begin{figure}[!b]
    \centering
    \vspace{-0.3cm}
    \scalebox{1.1}{\includegraphics[]{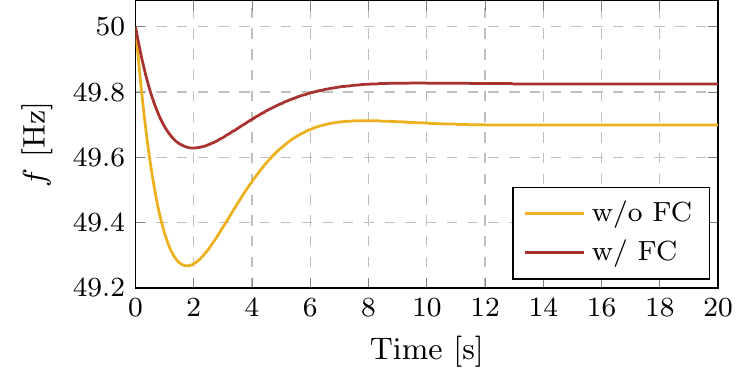}}
    \caption{Frequency evolution of the CoI with and without FC for 20 seconds after the fault instance at $t=67\,\mathrm{h}$.}
    \label{fig:deltaf_results}
\end{figure}

Some insightful conclusions can be drawn from Fig.~\ref{fig:IC_lim_values}, where the difference between the actual values of the frequency metrics and the respective ENTSO-E thresholds are depicted. For this purpose, we define a constraint gap $\eta$ as a measure of the relative constraint distance to its limit, e.g., $\eta_\mathrm{nadir}=\Delta f_\mathrm{max}/\Delta f_\mathrm{lim}-1$. After the completion of unit commitment, the constraints are re-evaluated using the obtained $F_g$, $R_g$ and $M$ values in order to determine which frequency criteria becomes binding at the instance of the fault. A negative constraint gap corresponds to the non-binding constraint, i.e., the specific frequency criteria being met. It should be noted that the positive values of $\eta$ for $t<67\,\mathrm{h}$ indicate that the frequency threshold would be violated if a fault occurs. No action is required however, considering that in this case study we assume that the contingency can only occur in hour $t=67$. The Fig.~\ref{fig:IC_lim_values} indicates that without explicit frequency metric constraints all of these constraints would be violated. Moreover, it can be observed that when including the frequency constraints, the RoCoF constraint is closest to its limit - corresponding to the smallest constraint gap - and thus binding. The constraint gap difference between the two scenarios at hour $t=67$ clearly highlights the importance of including the frequency constraints in UC in order to avoid large frequency excursions and undesired triggering of protection and UFLS schemes. The same observations are also reflected in Fig.~\ref{fig:deltaf_results} through time-domain frequency response of the system. Understandably, the values of RoCoF, nadir and steady-state deviation are reduced compared to the scenario without frequency constraints, such that all of the ENTSO-E criteria are fulfilled.

Finally, we investigate the economic impact of including the frequency constraints into the stochastic UC model. The breakdown of operational costs for day 3 is presented in Table~\ref{tab:price_breakdown}. The addition of frequency constraints leads to a $5\,\%$ increase in total expected system costs and a significant increase in start-up costs by $185\,\%$. This is due to six extra generators being turned on for providing inertia at the period of a potential generation failure, as shown in Table~\ref{tab:commitment_comp}. A large increase is also seen in reserve scheduling costs, as the reserves are now not only scheduled to cover wind power uncertainty but also for possible contingencies. The change in cost is of course highly dependent on the specific system and the considered contingencies. 

\section{Conclusion} \label{sec:Conclusion}

This paper includes frequency constraints in the UC problem of a system with large wind power penetration in order to investigate the impact of frequency dynamics on unit scheduling. By employing the analytic expressions for post-contingency frequency response of a multi-machine system, we define a set of constraints reflecting the frequency nadir, RoCoF and quasi steady-state deviation. The highly non-linear frequency nadir constraint is linearized using two approaches: (i) a PWL technique adapted from the literature; and (ii) a proposed simple and efficient method for extracting bounds on decision variables of interest, which is shown to be computationally superior to PWL. Using the latter approach, the stochastic UC problem is formulated as an MILP, with an objective of minimizing the expected system costs against wind power production and generation outage uncertainties.

\begin{table}[!t]
\renewcommand{\arraystretch}{1.2}
\caption{Unit commitment costs $\mathrm[\$]$ breakdown on day 3.}
\label{tab:price_breakdown}
\noindent
\centering
    \begin{minipage}{\linewidth} 
    \renewcommand\footnoterule{\vspace*{-5pt}} 
    \begin{center}
        \begin{tabular}{ c || c | c | c | c }
            \toprule
            \textbf{Scenario} & \textbf{Total costs} & \textbf{Start-up} & \textbf{Operation} & \textbf{Reserves} \\ 
            \cline{1-5}
            w/o FC & $410\,545$ & $788$ & $405\,561$ & $4\,196$\\
            \arrayrulecolor{black!30}\hline
            w/ FC & $432\,383$ & $2\,248$  & $420\,886$ & $9\,249$ \\
            \arrayrulecolor{black!30}\hline
            Difference  & $5.32\,\%$  & $185\,\%$  & $3.78\,\%$  & $120\,\%$  \\
            \arrayrulecolor{black}\bottomrule
        \end{tabular}
        \end{center}
    \end{minipage}
    \vspace{-0.3cm}
\end{table}

Our results show that the inclusion of frequency constraints in the UC model significantly affects the dispatch of synchronous generators and consequently the expected system costs. Indeed, during anticipated critical events such as the loss of generation, additional synchronous machines are needed for providing sufficient inertia and damping in the process of frequency containment. Such actions lead to a drastic increase in the UC costs, especially start-up and reserve scheduling, which poses a new challenge as the operator must find a way to remunerate the units committed for the sole purpose of frequency regulation. This is an exciting avenue for future work.

\appendices
\section{}  \label{appendixA}

\begin{figure}[!t]
    \centering
    \scalebox{0.9}{\includegraphics{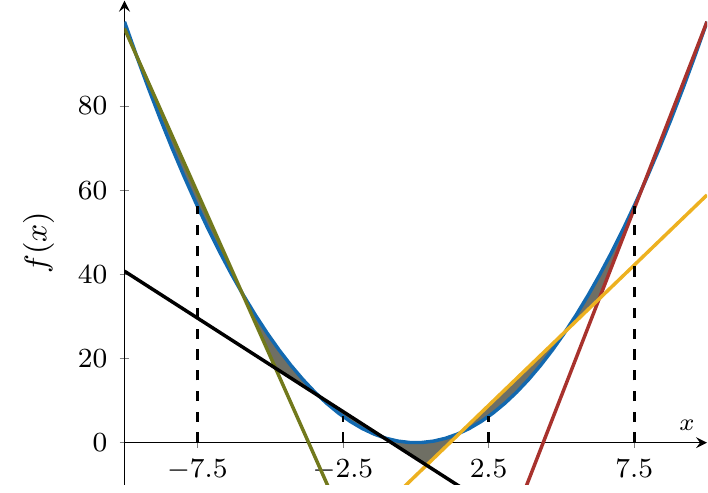}}
    \caption{Illustration of the PWL method on a 2-D function.}
    \label{fig:pwl_example}
\end{figure}

For clarity, we visually illustrate in Fig.~\ref{fig:pwl_example} how the PWL optimization problem from \eqref{eq:3d_pwl_obj} is solved. In this example there are four evaluation points set at -7.5, -2.5, 2.5 and 7.5, and the respective function is approximated with four segments. At each evaluation point the model identifies the segment closest to the original curve, and subsequently aims to minimize the overall shaded area. The proposed technique can then be expanded and employed on a function of three variables, as we have done with the expression for frequency nadir.

 References section
\bibliographystyle{IEEEtran}
\bibliography{bibliography}

\end{document}